\documentclass{article}

\usepackage{xcolor}
\usepackage{hyperref}
\usepackage{geometry,amssymb,bm,fullpage,amsmath,minipage-marginpar,authblk}
\usepackage{verbatim}
\usepackage{color,xcolor}
\usepackage{enumitem}
\usepackage{graphicx,appendix}

\usepackage{pgfplots}
\usepackage{tikz}
\usepgfplotslibrary{groupplots}
\usetikzlibrary{calc}
\pgfplotsset{
  compat=1.8,
  minor grid style={dashed,red},
  major grid style={dotted,green!50!black}
}

\newcommand{\Rmnum}[1]{\expandafter\@slowromancap\romannumeral #1@}
\makeatletter
 \newcommand{\cR}{\mathcal{R}}
 \newcommand{\sL}{\mathsf{\Phi}}
\newtheorem{assumption}{Assumption}[section]
\newtheorem{theorem}{Theorem}[section]

\newtheorem{example}{Example}[section]

\newtheorem{remark}{Remark}[section]

\newenvironment{pf}{{\noindent \it \bf Proof. }}{{\hfill$\Box$}\\}

\makeatother

\begin{document}

\title{Error Lower Bounds of Constant Step-size Stochastic Gradient Descent}
\author[1]{Zhiyan Ding}
\author[2]{Yiding Chen}
\author[1]{Qin Li}
\author[2]{Xiaojin Zhu}
\affil[1]{Department of Mathematics, University of Wisconsin-Madison}
\affil[2]{Department of Computer Sciences, University of Wisconsin-Madison}
\date{}    

\maketitle
%

%





\begin{abstract}
Stochastic Gradient Descent (SGD) plays a central role in modern machine learning. While there is extensive work on providing error upper bound for SGD, not much is known about SGD error lower bound. In this paper, we study the convergence of constant step-size SGD. We provide error lower bound of SGD for potentially non-convex objective functions with Lipschitz gradients. To our knowledge, this is the first analysis for SGD error lower bound without the strong convexity assumption. We use experiments to illustrate our theoretical results.
\end{abstract}

\section{Introduction}
Stochastic Gradient Descent (SGD) is one of the most popular optimization algorithms in modern machine learning
~\cite{bottou2010large,bottou2018optimization,zhang2004solving,badrinarayanan2017segnet,zinkevich2010parallelized}. SGD enjoys computational efficiency and is easy to implement.

In terms of theoretical guarantee, there are extensive studies on analyzing the error upper bound for SGD. 
~\cite{zhang2004solving} studies SGD on regularized forms of linear prediction methods, where an error upper bound is given. A non-asymptotic convergence analysis is given in~\cite{moulines2011non}.~\cite{nguyen2018sgd} studies the convergence of SGD with decaying step-size SGD.~\cite{rakhlin2011making} gives convergence analysis for averaged SGD.~\cite{nguyen2018tight} discusses the tightness of existing theoretic upper bound and provides a particular objective function whose SGD error almost touches this upper bound.~\cite{gower2019sgd} gives general convergence analysis under expected smooth assumption.
However, there is much less work on the SGD error lower bound.~\cite{dieuleveut2017bridging} provides the lower bound for averaged SGD with strongly convex objective functions.~\cite{jentzen2019lower} quantifies the convergence speed of SGD and gives error lower bound but they are considering SGD with decaying step-size. Their discussion is limited to quadratic objective functions and lack explicit formula of coefficient. When compared to prior work, our work is more general in the sense that we consider a family of objective function with Lipschitz gradient (both \textbf{convex} and \textbf{non-convex} functions are included). 

The error lower bound of SGD is of interest for several reasons.
On one hand, the error lower bound sharpens our understanding of SGD.~\cite{zhang2004solving} givens an error upper bound for constant step-size SGD when applied to a linear regression learning task. But there exists a term that does not vanish as iteration goes to infinity. This term depends on the step-size $\eta$ and will approach $0$ as $\eta \rightarrow 0$.
Our error lower bound accounts for the existence of such a term.
On the other hand, SGD error lower bound might be an interesting topic to the adversarial machine learning community. For training data poisoning attacks where an attacker changes training data before learning, a learner using constant step-size SGD will be partially immune to the attack in the sense that the attacker cannot precisely control the resulting model.

Our contributions include the following:
\begin{enumerate}
\item We give asymptotic and non-asymptotic analysis for SGD error lower and upper bounds when the objective function is strongly convex and has Lipschitz gradients (Theorem~\ref{Lemma:iteration}, \ref{ThCB1});
\item We give asymptotic analysis for SGD error lower bound when the objective function is potentially non-convex and has Lipschitz gradients (Theorem \ref{ThCB2}). We also give a loose but uniform constant lower bound (Theorem \ref{constantlowerbound}).
\end{enumerate}

\section{Preliminaries}

Stochastic gradient descent (SGD) is an iterative method for optimizing an objective function with the following structure
\begin{equation}\label{OMP}
\min_{x\in\mathbb{R}^d} f(\theta):=\frac{1}{J}\sum^J_{j=1}f_j(\theta).
\end{equation}
where $f, f_j:\mathbb{R}^d\rightarrow\mathbb{R}$ for $j=1,\dots,J$. In machine learning applications, $f$ is the total loss function whereas each $f_j$ represents the loss due to the $j$-th training sample. $\theta$ is a vector of trainable parameters and $J$ is the training sample size, which is typically very large. Throughout the paper we call $\theta^\dagger$ the optimal solution. Under the assumption that all objective functions are differentiable, we have
\begin{equation}\label{udagger}
\nabla f(\theta^\dagger) = 0\,.
\end{equation}

At each iteration, SGD randomly picks one of $\{f_j\}$ and performs gradient descent along the gradient provided by the chosen $f_j$. More explicitly, to obtain the $k+1$-th step solution, one has the following formula:
\begin{equation}\label{Istep}
\theta_{k+1}=\theta_k-\eta \nabla f_{\gamma_k}(\theta_k)\,,
\end{equation}
where $\eta>0$ is the step-size, and $\gamma_k$ is an $i.i.d.$ random variable from the uniform distribution on $\{1,2,\dots,J\}$.
In this paper we focus on SGD with constant step-size.

Comparing to Gradient Descent (GD), which requires gradients of the entire list of objective functions, SGD requires fewer number of gradients, and thus is computationally much cheaper per iteration. However, since stochasticity is involved, the convergence analysis is significantly more complicated. In particular, since $\nabla f_j (\theta^\dagger)$ may not be zero (despite $\nabla f(\theta^\dagger)=0$), SGD cannot terminate itself even at the global minimum. In the strongly convex case, we show that instead of converging to the minimum point, SGD forms a ``fuzzy error ball'' around it. 

We assume the gradients of $f$ and each $f_j$ to be Lipschitz:
\begin{assumption}\label{ass:smooth}
$\lambda$-smooth assumption:  there exist constants $\{\lambda_{\max,j}\}^{J}_{j=0}>0$ such that
\begin{equation}\label{eqn:assum_Lip}
    \begin{aligned}
|\nabla f(\theta_1)-\nabla f(\theta_2)| &\leq \lambda_{\max,0}|\theta_1-\theta_2| \\
|\nabla f_j(\theta_1)-\nabla f_j(\theta_2)| &\leq \lambda_{\max,j}|\theta_1-\theta_2|\;(\forall j)\,.
\end{aligned}
\end{equation}
\end{assumption}
In this paper, we use $| \cdot |$ to denote the $2$-norm of a vector.

Our analysis is divided into two situations: \textbf{Lipschitz} (Assumption~\ref{ass:smooth}) and \textbf{strongly convex} (Assumptions~\ref{ass:convex}); \textbf{Lipschitz} and potentially \textbf{non-convex}.
For the first situation, we will assume:
\begin{assumption}\label{ass:convex}
Strongly convex assumption:  there exist a constant $\lambda_{\min}>0$  such that
\begin{equation}\label{eqn:assum_conv}
\left\langle \nabla f(\theta_1)-\nabla f(\theta_2), \theta_1-\theta_2\right\rangle\geq \lambda_{\min}|\theta_1-\theta_2|^2\,,\\
\end{equation}
\end{assumption}
Under the strong convexity assumption the solution $\theta^\dagger$ satisfying~\eqref{udagger} is the unique solution to~\eqref{OMP}.

We introduce the following notations:
\[
D_0=\left(\frac{1}{J}\sum^J_{j=1}|\nabla f_j(\theta^\dagger)|^2\right)^{1/2}\,,\quad\Lambda=\left(\frac{1}{J}\sum^J_{j=1}\lambda^2_{\max,j}\right)^{1/2}\,.
\]

We now show two examples that satisfy both Assumption~\ref{ass:smooth} and~\ref{ass:convex}:
\begin{example}
\textbf{Linear Regression} Let $\{(\mathbf{x}_j, y_j)\}_{j = 1}^J$ be a training set, where $\mathbf{x}_j \in \mathbb{R}^{d \times 1}$ and $y_j \in \mathbb{R}$. 
We assume $J >> d$ and $X = (\mathbf{x}_1, \ldots, \mathbf{x}_J) \in \mathbb{R}^{d\times J}$ has full row rank.
We consider a linear regression problem. 
Let $f_j(\theta) = 1/2(\theta^\top\mathbf{x}_j - y_j)^2$.  
Then $f(\theta) = 1/(2J) \left|X^\top \theta - \mathbf{y}\right|^2$, where $\mathbf{y}  = (y_1, \ldots, y_J)^\top \in \mathbb{R}^{J \times 1}$.
In this case,~\eqref{OMP} has a closed-form solution and $\theta^\dagger$ is given by $(XX^\top)^{-1}X\mathbf{y}$. 
$\nabla f_j(\theta) = \mathbf{x}_j (\mathbf{x}_j^\top \theta - y_j)$, $\nabla f(\theta) =  X(X^\top \theta - \mathbf{y})/J$, $\nabla^2 f_j(\theta) = \mathbf{x}_j \mathbf{x}_j^\top$, $\nabla^2 f(\theta) =  XX^\top /J$. 
Thus we can choose $\lambda_{\max,0} = \sigma_{\max} (XX^\top)/J$, $\lambda_{\max,j} = \mathbf{x}_j^\top \mathbf{x}_j$ and $\lambda_{\min} = \sigma_{\min}(XX^\top)/J$. We use $\sigma_{\min}(\cdot)$ and $\sigma_{\max}(\cdot)$ to denote the smallest and largest eigenvalues.

\end{example}

\begin{example}
\textbf{Logistic Regression with $L_2$ Regularization} Let $\{(\mathbf{x}_j, y_j)\}_{j = 1}^J$ be a training set, where $\mathbf{x}_j \in \mathbb{R}^{d \times 1}$ and $y_j \in \{0,1\}$. 
The objective function of Logistic Regression with regularization (of weight 1) is
\begin{equation}
 \theta^\top \theta/2-1/J\sum_{j=1}^J(y_j \log(S(\theta^\top \mathbf{x}_j)) + (1-y_j) \log(1- S(\theta^\top \mathbf{x}_j))),
\end{equation}
where $S(\cdot)$ is the Sigmoid function defined by $S(x) = 1/(1+e^{-x})$.
In this case, 
\begin{equation}
f_j(\theta) = -(y_j \log(S(\theta^\top \mathbf{x}_j)) + (1-y_j) \log(1- S(\theta^\top \mathbf{x}_j))) + \theta^\top \theta/2,
\end{equation}
Then 
\begin{equation}
\nabla f_j(\theta)  =  (S(\theta^\top \mathbf{x}_j) - y_j) \mathbf{x}_j + \theta 
\end{equation}
By 
\begin{align*}
|\nabla f_j(\theta_1) - \nabla f_j(\theta_2)| &\le |S(\theta_1^\top \mathbf{x}_j) - S(\theta_2^\top \mathbf{x}_j)| |\mathbf{x}_j| + |\theta_1 - \theta_2| \\
&\le |\theta_1^\top \mathbf{x}_j - \theta_2^\top \mathbf{x}_j| |\mathbf{x}_j| + |\theta_1 - \theta_2|  \\
& \le (|\mathbf{x}_j|^2+1) |\theta_1 - \theta_2|,
\end{align*}
and
\begin{equation}
\nabla^2 f_j(\theta)  = S(\theta^\top \mathbf{x}_j) (1 - S(\theta^\top \mathbf{x}_j))\mathbf{x}_j \mathbf{x}_j ^\top + I,
\end{equation}
 we can choose 
\begin{equation*}
\lambda_{\max,j} = (|\mathbf{x}_j|^2 + 1), \lambda_{\max,0} = 1+1/J\sum_{j=1}^J|\mathbf{x}_j|^2, \lambda_{\min} = 1
\end{equation*}
\end{example}




\section{Main Results}
We present our main theory in this section. It is our goal to show that SGD is bounded away from the minimum.
The quantity of interest is: 
\begin{equation}\label{Def:WK}
\cR_k=\mathbb{E}\left|\theta_k-\theta^\dagger\right|^2\,.
\end{equation}
It is the expected squared distance between the $k$-iteration solution and the target $\theta^\dagger$. The expectation takes into account the randomness in SGD.


The analysis builds upon two steps: in Step 1, we derive the upper and lower bound of the iterative formula that updates $\cR_{k+1}$ from $\cR_k$. This step heavily depends on rewriting the SGD formula. We present the result in Subsection \ref{Sec:recuresive}. In Step 2, we analyze the long time behavior of the iterative formula. The behavior depends on the strong convexity of the cost function, and thus we separate the discussion for the strongly convex and non-convex case, and present them in Subsection \ref{Sec:convex} and \ref{Sec:nonconvex}, respectively.

\subsection{Key Recursive Inequalities}\label{Sec:recuresive}
The derivation of the updating formula for $\cR_k$ comes directly from SGD~\eqref{Istep}, which we rewrite into:
\begin{equation}\label{eqn:lstep_V}
\theta_{k+1}=\theta_k-\eta\nabla f(\theta_k)+\eta V\left(\theta_k\right)\,,
\end{equation}
where the variation
\[
V\left(\theta\right)=\nabla f(\theta)-\nabla f_{\gamma_k}(\theta)\,.
\]
Due to the randomness involved in the selection of $\gamma_k$, $V$ is also a random variable with explicitly expressible mean and variance:
\[
\begin{aligned}
&\mathbb{E}\left(V\left(\theta\right)\right)=\overrightarrow{0}\,,\\ &\mathrm{Var}\left(V\left(\theta\right)\right)=\frac{1}{J}\sum^{J}_{j=1}\left(\nabla f(\theta)-\nabla f_j(\theta)\right)\left(\nabla f(\theta)-\nabla f_j(\theta)\right)^\top\,.
\end{aligned}
\]

Building upon these, we have the following theorem on the recursive formula.
\begin{theorem}\label{Lemma:iteration}
The following upper and lower bounds hold true if SGD uses fixed step-size $\eta<\frac{1}{\lambda_{\max,0}}$:
\begin{itemize}
\item If the objective functions are $\lambda$-smooth (satisfy Assumption~\ref{ass:smooth}), then
\begin{equation}\label{infinequality}
\cR_{k+1}-\cR_k\geq \left[-2\eta\lambda_{\max,0}\right]\cR_k-2\eta^2\Lambda D_0\cR^{1/2}_{k}+\eta^2D^2_0\,.
\end{equation}
\item In addition, if the objective functions are also strongly convex (satisfy both Assumption~\ref{ass:smooth} and~\ref{ass:convex}), then:
\begin{align}\label{supinequality}
\cR_{k+1}-\cR_k\leq & \left[-2\eta\lambda_{\min}+\eta^2\left(\lambda^2_{\max,0}+\Lambda^2-\lambda^2_{\min}\right)\right]\cR_k +2\eta^2\Lambda D_0\cR^{1/2}_{k}+\eta^2D^2_0\,.
\end{align}
\end{itemize}
\end{theorem}

\begin{remark}
There are several comments:
\begin{itemize}
\item Equation~\eqref{infinequality} gives the lower bound for $\cR_{k+1}$ while~\eqref{supinequality} gives the upper bound. The lower bound relies only on the $\lambda$-smooth assumption while the upper bound furthermore requires strong convexity. This is in line with the standard GD analysis.
\item For very small $\eta$, the $\eta^2$ terms in~\eqref{infinequality} and~\eqref{supinequality} become negligible, and the linear term on $\cR_k$ dominates the estimates.
\item In the recursive formula for the classical GD analysis the new-step error linearly depends on the old-step error. We do have an extra square root term $\sqrt{\cR_k}$. This is exactly because the randomness adds a direction different from the gradient descent direction. To handle this new direction, the H\"older inequality is utilized and the power of $\cR_k$ is thus sacrificed.
\item In the second statement of the theorem, strongly convex property is imposed. This may not be necessary, see~\cite{gower2019sgd}. However, it is our main goal to understand the lower but not the upper bound, and thus we do not pursue a tighter condition in this paper.
\item The statement of the theorem still holds true for varying time step, assuming that $\eta_k$ satisfies the condition at every step: $\eta_k<\frac{1}{\lambda_{\text{max},0}}$ for all $k$.

\item For $J = 1$, SGD degenerates to GD and $D_0 = 0$. \eqref{supinequality} recover the classical analysis for GD.
\end{itemize}
\end{remark}

\begin{pf}
We start from~\eqref{eqn:lstep_V}, subtract $\theta^\dagger$ and take expectation of the $L_2$ norm and get:
\begin{equation}\label{Secondcal}
\begin{aligned}
\cR_{k+1} &= \mathbb{E}\left|\theta_{k+1}-\theta^\dagger\right|^2=\mathbb{E}\left|\theta_k-\theta^\dagger-\eta\nabla f(\theta_k)\right|^2+\eta^2\mathbb{E}\left|V(\theta_k)\right|^2\\
&=\mathbb{E}\left|\theta_k-\theta^\dagger-\eta\left(\nabla f(\theta_k)-\nabla f(\theta^\dagger)\right)\right|^2+\eta^2\mathbb{E}\left|V(\theta_k)\right|^2\\
&= \mathbb{E}\left|\theta_k-\theta^\dagger-\eta\left(\nabla f(\theta_k)-\nabla f(\theta^\dagger)\right)\right|^2+\frac{1}{J}\sum^J_{j=1}\mathbb{E}|\nabla f_j(\theta_k)|^2-\mathbb{E}|\nabla f(\theta_k)|^2\,.
\end{aligned}
\end{equation}
In the first equation we used the fact that $\mathbb{E}\left(V\left(\theta\right)\right)=0$, and in the second we use the fact that $\nabla f(\theta^\dagger)=0$. Lastly we expand $\mathbb{E}|V(\theta_k)|^2$ using the definition of $V$.

To obtain~\eqref{infinequality} and~\eqref{supinequality} amounts to give upper and lower bounds for the three terms on the right hand side.

To lower bound the first term, we note that using the triangle inequality:
\begin{equation}\label{eqn:term1_lower}
\begin{aligned}
\mathbb{E}\left|\theta_k-\theta^\dagger-\eta\left(\nabla f(\theta_k)-\nabla f(\theta^\dagger)\right)\right|^2\geq&\mathbb{E}\left(\left|\theta_k-\theta^\dagger\right|-\eta\left|\nabla f(\theta_k)-\nabla f(\theta^\dagger)\right|\right)^2\\
\geq&\mathbb{E}\left(\left(1-\eta\lambda_{\max,0}\right)\left|\theta_k-\theta^\dagger\right|\right)^2\\
\geq&\left(1-\eta\lambda_{\max,0}\right)^2\mathbb{E}\left|\theta_k-\theta^\dagger\right|^2 = (1-\eta\lambda_{\max,0})^2\cR_k\,,
\end{aligned}
\end{equation}
where the second inequality comes from the $\lambda$-smooth assumption and the requirement on $\eta$.

To upper bound the first term, we note that:
\begin{equation*}
\begin{aligned}
&\left|\theta_k-\theta^\dagger-\eta\left(\nabla f(\theta_k)-\nabla f(\theta^\dagger)\right)\right|^2\\
=&\langle \theta_k-\theta^\dagger\,, \theta_k-\theta^\dagger\rangle - 2\eta \langle \theta_k-\theta^\dagger\,, \nabla f(\theta_k)-\nabla f(\theta^\dagger)\rangle  +\eta^2\langle \nabla f(\theta_k)-\nabla f(\theta^\dagger)\,,\nabla f(\theta_k)-\nabla f(\theta^\dagger)\rangle\,,
\end{aligned}
\end{equation*}
which gives
\begin{equation}\label{eqn:term1_upper}
\begin{aligned}
\mathbb{E}\left|\theta_k-\theta^\dagger-\eta\left(\nabla f(\theta_k)-\nabla f(\theta^\dagger)\right)\right|^2\leq&\left(1-2\eta\lambda_{\min}+\eta^2\lambda^2_{\max,0}\right)\mathbb{E}\left|\theta_k-\theta^\dagger\right|^2=\left(1-2\eta\lambda_{\min}+\eta^2\lambda^2_{\max,0}\right)\cR_k\,.
\end{aligned}
\end{equation}
by the strong convexity condition.

To estimate the third term in \eqref{Secondcal}, we note that on one hand, by using the $\lambda$-smoothness assumption \ref{ass:smooth}, we have:
\begin{equation}\label{eqn:term2_2}
\begin{aligned}
\mathbb{E}|\nabla f(\theta_k)|^2&=\mathbb{E}|\nabla f(\theta_k)-\nabla f(\theta^\dagger)|^2\leq \lambda^2_{\max,0}\mathbb{E}\left|\theta_k-\theta^\dagger\right|^2 =  \lambda^2_{\max,0}\cR_k\,,
\end{aligned}
\end{equation}
and on the other, by using the strong convexity property \eqref{eqn:assum_conv}, we have:
\begin{equation}\label{eqn:term2_3}
\begin{aligned}
\mathbb{E}|\nabla f(\theta_k)|^2&=\mathbb{E}|\nabla f(\theta_k)-\nabla f(\theta^\dagger)|^2\geq \lambda^2_{\min}\mathbb{E}\left|\theta_k-\theta^\dagger\right|^2= \lambda^2_{\min}\cR_k\,.
\end{aligned}
\end{equation}
To estimate the second term in~\eqref{Secondcal}, we insert $\theta^\dagger$ to have: 
\begin{equation}\label{eqn:term2_1_upper}
\begin{aligned}
&\frac{1}{J}\sum^J_{j=1}\mathbb{E}|\nabla f_j(\theta_k)|^2 = \frac{1}{J}\sum^J_{j=1}\mathbb{E}|\nabla f_j(\theta_k)-\nabla f_j(\theta^\dagger)+\nabla f_j(\theta^\dagger)|^2\\
=&\frac{1}{J}\sum^J_{j=1}\mathbb{E}|\nabla f_j(\theta^\dagger)|^2+\frac{1}{J}\sum^J_{j=1}\mathbb{E}|\nabla f_j(\theta_k)-\nabla f_j(\theta^\dagger)|^2+\frac{2}{J}\sum^J_{j=1}\mathbb{E}\left\langle \nabla f_j(\theta^\dagger),\nabla f_j(\theta_k)-\nabla f_j(\theta^\dagger)\right\rangle\,.
\end{aligned}
\end{equation}
Noting that $\frac{1}{J}\sum_{j=1}^J\mathbb{E}|\nabla f_j(\theta^\dagger)|^2 = D^2_0$,
\[
0\leq \frac{1}{J}\sum^J_{j=1}\mathbb{E}|\nabla f_j(\theta_k)-\nabla f_j(\theta^\dagger)|^2\leq \Lambda^2\mathbb{E}|\theta_k-\theta^\dagger|^2
\]
and that
\begin{equation}
\begin{aligned}
\frac{1}{J}\sum^J_{j=1}\mathbb{E}\left\langle \nabla f_j(\theta^\dagger),\nabla f_j(\theta_k)-\nabla f_j(\theta^\dagger)\right\rangle\leq& \frac{1}{J}\sum^J_{j=1}\left|\nabla f_j(\theta^\dagger)\right|\mathbb{E}\left|\nabla f_j(\theta_k)-\nabla f_j(\theta^\dagger)\right|\\
\leq& \frac{1}{J}\sum^J_{j=1}|\nabla f_j(\theta^\dagger)|\mathbb{E}\left(\lambda_{\max,j}\left|\theta_k-\theta^\dagger\right|\right)\\
\leq& \frac{1}{J}\sum^J_{j=1}|\nabla f_j(\theta^\dagger)|\lambda_{\max,j} \left(\mathbb{E}|\theta_k-\theta^\dagger|^2\right)^{1/2}\\
\leq& \left(\frac{1}{J}\sum_{j=1}^J|\nabla f_j(\theta^\dagger)|^2\right)^{1/2}\left(\frac{1}{J}\sum_{j=1}^J\lambda^2_{\max,j}\right)^{1/2} \sqrt{\cR_k}\\
=&D_0{\Lambda}\sqrt{\cR_k}\,.
\end{aligned}
\end{equation}
The second to last inequality relies on the $\lambda$-smooth assumption and the last one comes from Cauchy-Schwartz inequality. Putting together these estimates leads to the upper and lower bound of the second term in~\eqref{Secondcal}:
\begin{equation}\label{eqn:term2}
D^2_0- 2\Lambda D_0\sqrt{\cR_k}\leq \frac{1}{J}\sum^J_{j=1}\mathbb{E}|\nabla f_j(\theta_k)|^2\leq D^2_0 + \Lambda^2\cR_k + 2\Lambda D_0\sqrt{\cR_k}\,.
\end{equation}

To finish the proof, we combine~\eqref{eqn:term1_lower}, \eqref{eqn:term2} and~\eqref{eqn:term2_2} for the lower bound:
\[
\cR_{k+1}\geq \left[1-2\eta\lambda_{\max,0}\right]\cR_k-2\eta^2\Lambda D_0\cR^{1/2}_{k}+\eta^2D^2_0\,.
\]
With the strongly convex assumption \ref{ass:convex}, combining~\eqref{eqn:term2},~\eqref{eqn:term2_3} and~\eqref{eqn:term1_upper}, we have:
\[
\begin{aligned}
\cR_{k+1}\leq&\left[1-2\eta\lambda_{\min}+\eta^2\left(\lambda^2_{\max,0}+\Lambda^2-\lambda^2_{\min}\right)\right]\cR_k+2\eta^2\Lambda D_0\cR^{1/2}_{k}+\eta^2D^2_0\,,
\end{aligned}
\]
concluding the proof.
\end{pf}


\subsection{Strongly convex Objective Function}\label{Sec:convex}
If the objective function is strongly convex, it is a well-known result that SGD converges to the optimal solution with decaying step-size (\cite{nguyen2018sgd,gower2019sgd,zhang2004solving,nguyen2018tight}). With fixed step-size, however, it is widely believed that $\theta_k$ oscillates around the true solution with a noise determined by the step-size. We give a tight trajectory of $\cR_k$ in the following theorem based on the recursive formula above.
\begin{theorem}\label{ThCB1}
Suppose the objective functions in~\eqref{OMP} are $\lambda_{max,0}$-smooth and strongly convex (satisfying~\eqref{eqn:assum_Lip} and~\eqref{eqn:assum_conv}), then using SGD with step-size:
\begin{equation}\label{lambdacondition}
\eta<\frac{\lambda_{\min}}{\Lambda^2+\lambda^2_{\max,0}}\,,
\end{equation}
we have constants $C_{0-4}$ depending on $f$, $f_i$ and $\cR_0$ only so that:
\begin{enumerate}[label=$\bullet$,topsep=0pt]
\item \textbf{Non-asymptotic rate for all $k$}
\begin{equation}\label{supiteration}
\begin{aligned}
&\cR_{k}\leq \left[1-2\eta\lambda_{\min}+\eta^2\mathsf{\Phi}^2\right]^k\cR_0+C_2\eta\,,\\
&\cR_{k}\geq \left[1-2\eta\lambda_{\max}\right]^k\cR_0+C_3\eta\,
\end{aligned}
\end{equation}
\item \textbf{Non-asymptotic rate for large $k$:} Let $K_0=\left\lceil\frac{\log(\eta)}{\log(1-\eta\lambda_{\min})}\right\rceil$, then we have
\[
\cR_{K_0}\leq C_4\eta\,,
\]
and for $k\geq K_0$:
\begin{equation}\label{supiteration2}
\begin{aligned}
\cR_{k} \leq &\left[1-2\eta\lambda_{\min}+\eta^2\mathsf{\Phi}^2\right]^{k-K_0}\cR_{K_0}+\frac{\eta D^2_0+2C^{1/2}_4\eta^{3/2}\Lambda D_0}{2\lambda_{\min}-\eta\mathsf{\Phi}^2}\,,\\
\cR_k \geq &\left[1-2\eta\lambda_{\max,0}\right]^{k-K_0}\cR_{K_0}+\frac{\eta D^2_0-2C^{1/2}_4\eta^{3/2}\Lambda D_0}{2\lambda_{\max,0}}\,.
\end{aligned}
\end{equation}

\item \textbf{Asymptotic rate for large $k$:}
\begin{align}
&\limsup_{k\rightarrow\infty}\cR_{k}\leq \frac{\eta D^2_0+2C^{1/2}_4\eta^{3/2}\Lambda D_0}{2\lambda_{\min}-\eta\mathsf{\Phi}^2}\,,\label{supb}\\
&\liminf_{k\rightarrow\infty}\cR_{k}\geq \max\left\{\frac{\eta D^2_0-2C^{1/2}_4\eta^{3/2}\Lambda D_0}{2\lambda_{\max,0}},0\right\}\,.\label{infb}
\end{align}

\end{enumerate}
In the estimate,
\[
\mathsf{\Phi}=\sqrt{\Lambda^2+\lambda^2_{\max,0}-\lambda^2_{\min}}\,,
\]
and the constants $C_{0-4}$ can be made explicit:
\begin{align*}
& C_0(\eta) =\frac{\eta\Lambda D_0+\sqrt{\eta^2\Lambda^2D^2_0+\eta D^2_0\left(2\lambda_{\min}-\eta\mathsf{\Phi}^2\right)}}{2\lambda_{\min}-\eta\mathsf{\Phi}^2}\,, \\
& C_1(\eta)=\max\left\{\cR_0\,, C_0^2\right\}\,,\\
&C_2(\eta)=\frac{2\Lambda D_0\left[C_1+\left(D^2_0+1/2\right)\eta^2\right]^{1/2}+D^2_0}{2\lambda_{\min}-\eta\mathsf{\Phi}^2}\,, \\
&C_3(\eta)=\frac{-2\Lambda D_0\left[C_1+\left(D^2_0+1/2\right)\eta^2\right]^{1/2}+D^2_0}{2\lambda_{\max,0}}\,, \\
& C_4(\eta)=\cR_0+C_2(\eta)\,.
\end{align*}
\end{theorem}

\begin{remark}
Several comments are in order.
\begin{itemize}
\item If $\cR_0\sim O(1)$, by definition, except $C_0\sim O(\sqrt{\eta})$, all other constants $C_{1-4}$ are of $O(1)$ on $\eta$.

\item The asymptotic rate \eqref{supb},\eqref{infb} are not very tight when time step is not small enough. If $\eta$ further satisfies
\begin{equation*}
\frac{\eta^2\Lambda D_0}{1-2\eta\lambda_{\max,0}}\leq \frac{-2\eta\Lambda D_0+\sqrt{4\eta^2\Lambda^2 D^2_0+8\eta\lambda_{\max,0}D^2_0}}{4\lambda_{\max,0}}\,,
\end{equation*}
we could improve them as
\begin{align}
\limsup_{k\rightarrow\infty}\cR_{k}\leq &\left[\frac{2\eta\Lambda D_0+\sqrt{4\eta^2(\Lambda^2 D^2_0-\mathsf{\Phi}^2D^2_0)+8\eta D^2_0\lambda_{\min}}}{4\lambda_{\min}-2\eta\mathsf{\Phi}^2}\right]^2\,,\label{supbre}\\
\liminf_{k\rightarrow\infty}\cR_{k}\geq &\left[\frac{-2\eta\Lambda D_0+\sqrt{4\eta^2\Lambda^2 D^2_0+8\eta\lambda_{\max,0}D^2_0}}{4\lambda_{\max,0}}\right]^2\,.\label{infbre}
\end{align}
We proved \eqref{infbre} in the non-convex case \eqref{infb2}. The upper bound \eqref{supbre} can be proved using similar techniques, which we omitted.
\end{itemize}
\end{remark}

In Theorem \ref{ThCB1}, we use estimation of $\cR^{1/2}_k$. This will give us some higher order (w.r.t $\eta$) error in the trajectory track \eqref{supiteration},\eqref{supiteration2}. In real application, if time step $\eta$ is not small enough, we can choose to calculate \eqref{infinequality},\eqref{supinequality} directly to track the trajectory of $\cR_k$.

It's clear from equation~\eqref{supiteration} and~\eqref{supiteration2} the impact of the error introduced due to the random initial guess is eliminated at an exponential rate. As $k$ increases, the second term dominates. It is a term of $\mathcal{O}(\eta)$, with the coefficient mainly determined by $D_0$, that reflects the influence of the stochasticity in the algorithm. This is a term that cannot be eliminated. Furthermore, comparing~\eqref{supb} with~\eqref{infb}, we see that in the leading order, for small $\eta$, approximately
\begin{align*}
\limsup_{k\rightarrow\infty}\cR_{k}\lesssim \frac{\eta D_0}{2\lambda_{\min}}\,,\quad\text{while}\quad \liminf_{k\rightarrow\infty}\cR_{k}\gtrsim \frac{\eta D_0}{2\lambda_{\max,0}}\,.
\end{align*}
This implies in the strongly convex case, if we run SGD forever, the expectation of the error will only be determined by the total loss function $f$ and the step-size $\eta$.

\subsection{Potentially Non-Convex Objective}\label{Sec:nonconvex}

Without the strongly convex assumption \ref{ass:convex}, it's hard to track the trajectory of $\cR$. The main reason is we only have one direct inequality \eqref{infinequality} in Theorem \ref{Lemma:iteration}. However, the lower bound of $\liminf$ has a relaxed requirement on strong convexity and the properties of objective function, and thus a quantitative result is still available for the lower bound of $\liminf$. The following theorem gives us the asymptotic lower bound of $\cR$:
\begin{theorem}\label{ThCB2}
Suppose the objective functions are $\lambda_{\max,0}$-smooth (satisfying~\eqref{eqn:assum_Lip}), then using SGD with step-size $\eta$ such that:
\[
\eta<\frac{1}{2\lambda_{\max,0}},
\]
there exists $z_0$ so that
\begin{enumerate}[label=$\bullet$,topsep=0pt]
\item \textbf{Asymptotic rate of $\limsup$:}
\begin{equation}\label{supb2}
\limsup_{k\rightarrow\infty}\cR_{k}\geq z_0^2\,.
\end{equation}

\item  \textbf{Asymptotic rate of $\liminf$:} If $\eta$ is small enough such that
\begin{equation}\label{smallcondition2}
\frac{\eta^2\Lambda D_0}{1-2\eta\lambda_{\max,0}}\leq z_0\,,
\end{equation}
then
\begin{equation}\label{infb2}
\liminf_{k\rightarrow\infty}\cR_{k}\geq z^2_0\,.
\end{equation}
\end{enumerate}
$z_0$ can be made explicit:
\begin{equation}\label{eqn:C_0C_1}
z_0(\eta)=\frac{-2\eta\Lambda D_0+\sqrt{4\eta^2\Lambda^2 D^2_0+8\eta\lambda_{\max,0}D^2_0}}{4\lambda_{\max,0}}.
\end{equation}
\end{theorem}
\begin{remark}
Several comments are in order:
\begin{itemize}
\item In this theorem, if we further have $\eta<\frac{\lambda_{\max,0}}{2\Lambda^2}$, then 
we have
\[
z_0\geq \sqrt{\frac{\eta D^2_0-\eta^{3/2}D_0\sqrt{2D^2_0\Lambda^2/\lambda_{\max,0}}}{2\lambda_{\max,0}}}.
\]
Therefore, when $\eta$ is small enough, the leading term of $\liminf$ \eqref{infb2} is the same as the strongly convex case \eqref{infb}.

\item Condition~\eqref{smallcondition2} is a very weak requirement on $\eta$. Indeed the left hand is a second order term of $\eta$ while the right hand side is $z_0$ which is about $\mathcal{O}(\eta^{1/2})$. This condition always holds true for small $\eta$.
\end{itemize}
\end{remark}

Compare Theorem \ref{ThCB2} with Theorem \ref{ThCB1}, the major difference is we lost the non-asymptotic estimation in the non-convex case. Without further assumption, it's hard to find upper bounds in each iteration like \eqref{supinequality}. This makes controlling negative term $-\cR^{1/2}_k$ and tracking tight lower bound difficult. However, using structure of the left hand side in \eqref{infinequality}, we can prove $z^2_0$ is also a loose lower bound for every iteration as the following theorem states:

\begin{theorem}\label{constantlowerbound}
If \eqref{smallcondition2} is true, then $\cR_k$ satisfies the following properties:
\begin{enumerate}[label=$\bullet$,topsep=0pt]
\item For any $k\geq 0$,
\[
\cR_{k}>\min\left\{\cR_{0},z^2_0\right\}
\]
\item For any $k\geq 0$, if $\cR_k<z^2_0$, then
\begin{equation}\label{alwayslarge1}
\cR_{k+1}>\cR_{k}
\end{equation}

\item If there exists $k^*$ such that $\cR_{k^*}\geq z^2_0$, then we have
\begin{equation}\label{alwayslarge2}
\cR_k\geq z^2_0,\quad \forall k>k^*\,.
\end{equation}

\end{enumerate}
\end{theorem}

In Theorem \ref{constantlowerbound}, equation~\eqref{alwayslarge1},\eqref{alwayslarge2} show if there is one step whose error, in expectation sense, is small, then the error will grow and if there is one step the error is above a certain threshold, in expectation sense, the error will be constantly above the threshold, leading to no-convergent case. 

In real applications, most likely $\cR_0>z^2_0$. By \eqref{alwayslarge2}, although we cannot get a tight lower bound of $\cR_k$ for any $k$, at least we know it will not decrease below the threshold $z^2_0$ in any steps.

\section{Experiments}
We now illustrate our theoretical results with experiments. 
In Section~\ref{sec:LR}, we use SGD to learn a linear regression model. 
In Section~\ref{sec:NonConvexExper}, we use SGD to minimize a $4$th order piecewise polynomial on $\mathbb{R}$. 

\subsection{SGD for Strongly Convex Quadratic Objective Function \label{sec:LR}}
In this experiment, we run SGD to learn a linear regression model $y = \theta^\top x + \epsilon$, where $\theta \in \mathbb{R}^{2\times 1}$. We will use this experiment to show the existence of an error lower bound and visualize the ``fuzzy ball'' around the global minimum.

The training data $\{(\mathbf{x}_j, y_j)\}_{j=1}^J \subset \mathbb{R}^{d \times 1} \times \mathbb{R}$ is generated as follows: the size of training set is $J = 30$, the dimension of the feature vector is $d = 2$. Each $\mathbf{x}_j$ is a random vector on the unit sphere generated as follows: 
first draw $\bm{x}_j \sim N(\bm{0}, I_{2 \times 2})$,
then normalize it into $\mathbf{x}_j = \bm{x}_j / |\bm{x}_j|$. 
We use $\theta^*= (-1.27,-0.49)^\top$ to generate the $y_j$'s.
For $j = 1, \ldots, J$, $y_j = \theta^{* \top}\mathbf{x}_j + \epsilon_j$, where $\epsilon_j \sim N(0, 0.1^2)$. $\theta^\dagger = (-1.27,-0.48)^\top$ is the Ordinary Least Square solution of this linear regression problem. 

\begin{figure}[ht]
\centering
\includegraphics[width=0.6\textwidth]{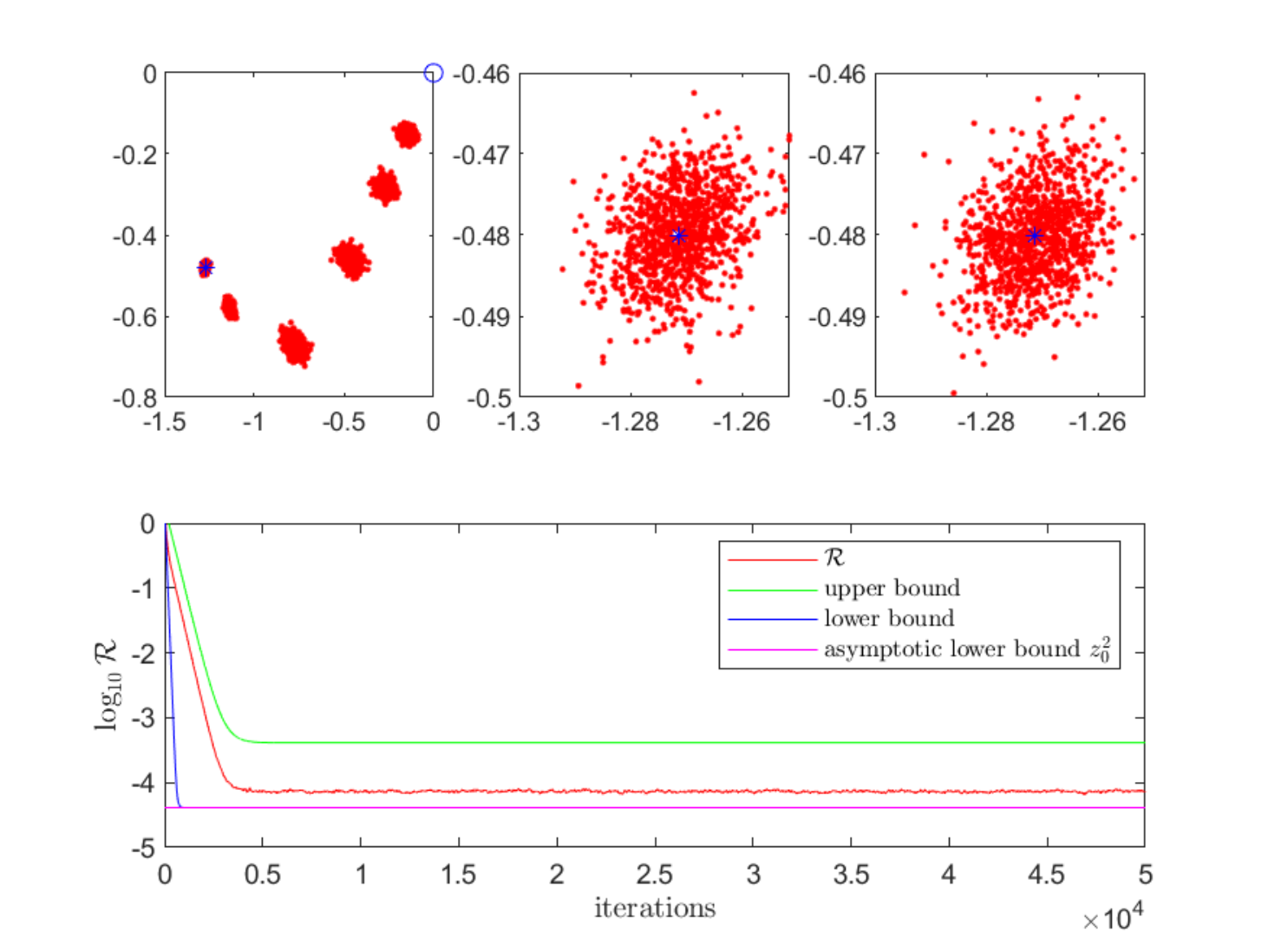}
\caption{\textbf{Upper panels}: upper left panel shows $\theta_k$ for $k = 25, 50, 100, 250, 1000, 5000, 10000, 50000$ in each trial; upper middle panel shows $\theta_{10000}$ in each trial; upper right panel shows $\theta_{50000}$ in each trial. In these panels, the blue circle is the initial point ($\theta_0$), the blue star is the global optimal ($\theta^\dagger$). We can only find $6$ clusters in the upper left panel because the clusters in the $5000, 10000$ and $50000$th iterations overlap. \textbf{Lower panel}: the SGD error curve for a strongly convex objective function\label{fig:exper1}.}
\end{figure}

When running SGD, we set the step-size to be $\eta = 0.01$ and set the maximum iteration to be $50000$. The initial point is $(0,0)^\top$. For the same data set, we run $1000$ trials. The only difference between these trials is due to the randomness in SGD. The error $\cR_k$ is an expectation, we estimate $\cR_k$ by averaging the error $|\theta_k - \theta^\dagger|^2$ over the $1000$ trials. We visualize the $\theta_k$'s for some particular iterations and plot the error curve in Figure~\ref{fig:exper1}. In the error curve, the lower and upper bounds are obtained by computing~\eqref{infinequality} and~\eqref{supinequality} in Theorem~\ref{Lemma:iteration}, the asymptotic lower bound $z_0^2$ is the result in Theorem~\ref{ThCB2}.

The upper panels in Figure~\ref{fig:exper1} shows the distribution of $\theta_k$ in particular iterations. Intuitively, the variance of $\theta_k$ does not decrease to zero as $t$ becomes large even after the SGD has stabilized.
The error curve in Figure~\ref{fig:exper1} shows a positive constant lower bound of the error, which demonstrates that the error of solution $\cR_k$ will not decrease to zero. Instead, the solution will be bounded away from a small fuzzy ball with a fixed radius in the sense of expected $L_2$ norm. So our results do not mean SGD will never touch the minimum point $\theta^\dagger$.

\subsection{SGD for Non-Convex Objective Function \label{sec:NonConvexExper}}
In this experiment, we run SGD to minimize a non-convex function. We use this experiments to show the existence of an error lower bound when the objective function is non-convex. 

\begin{figure}[ht]
\centering
\includegraphics[width=0.8\textwidth]{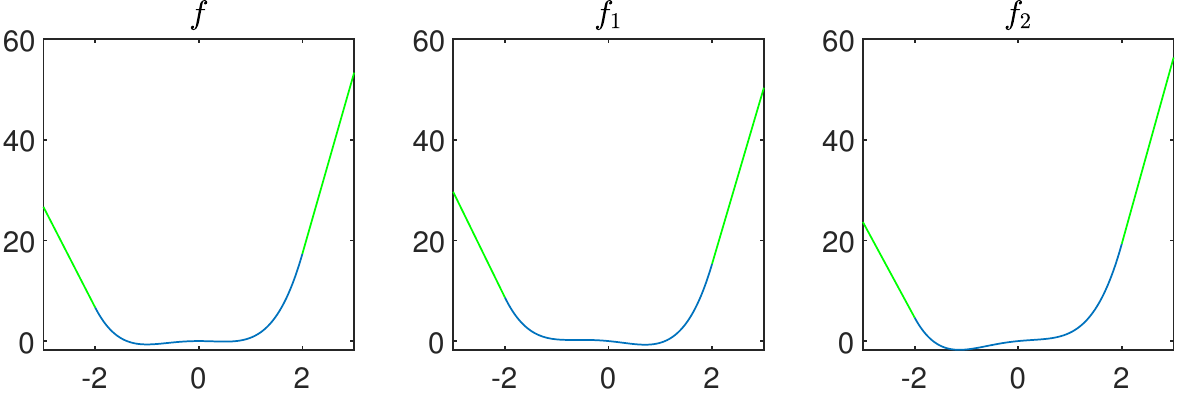}
\caption{$f,f_1,f_2$ for $\theta\in[-3, 3]$\label{fig:fplot}}
\end{figure}

For simplicity, we let the objective function be the following function: $f(\theta) = \theta^4 + 2/3 \theta^3 - \theta^2$, $f_1(\theta) = \theta^4 + 2/3 \theta^3 -\theta^2 +\theta$, $f_2(\theta) = \theta^4 + 2/3 \theta^3 -\theta^2 -\theta$ for $\theta \in [-2,2]$. To ensure that $f,f_1,f_2$ has Lipschitz derivatives, we define $f,f_1,f_2$ to be linear functions for $\theta \notin [-2,2]$ such that $f,f_1,f_2,f',f_1',f_2'$ are all continuous on $\mathbb{R}$. See Figure~\ref{fig:fplot}. It is clear that $f$ has two local minima at $\theta = -1$ and $\theta = 0.5$. $\theta = -1$ is the global minimum, i.e. $\theta^\dagger = -1$.

\begin{figure}[ht]
\begin{minipage}{.47\textwidth}
\centering
\includegraphics[width=1\textwidth]{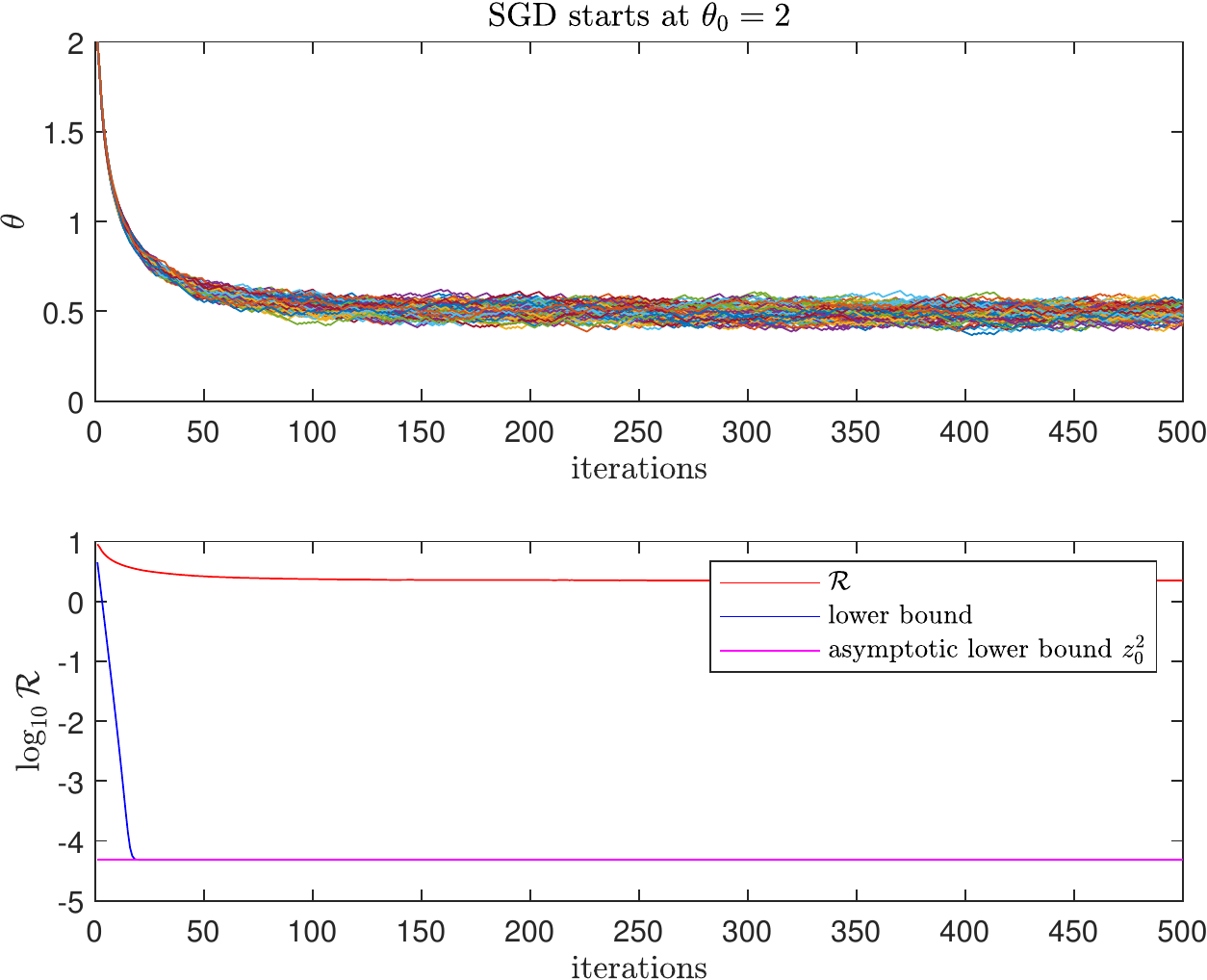}
\caption{$\theta$ trajectories and error curve for SGD starts at $\theta_0 = 2$, $500$ trials. \label{fig:err2_1}}
\end{minipage}
\begin{minipage}{.47\textwidth}
\centering
\includegraphics[width=1\textwidth]{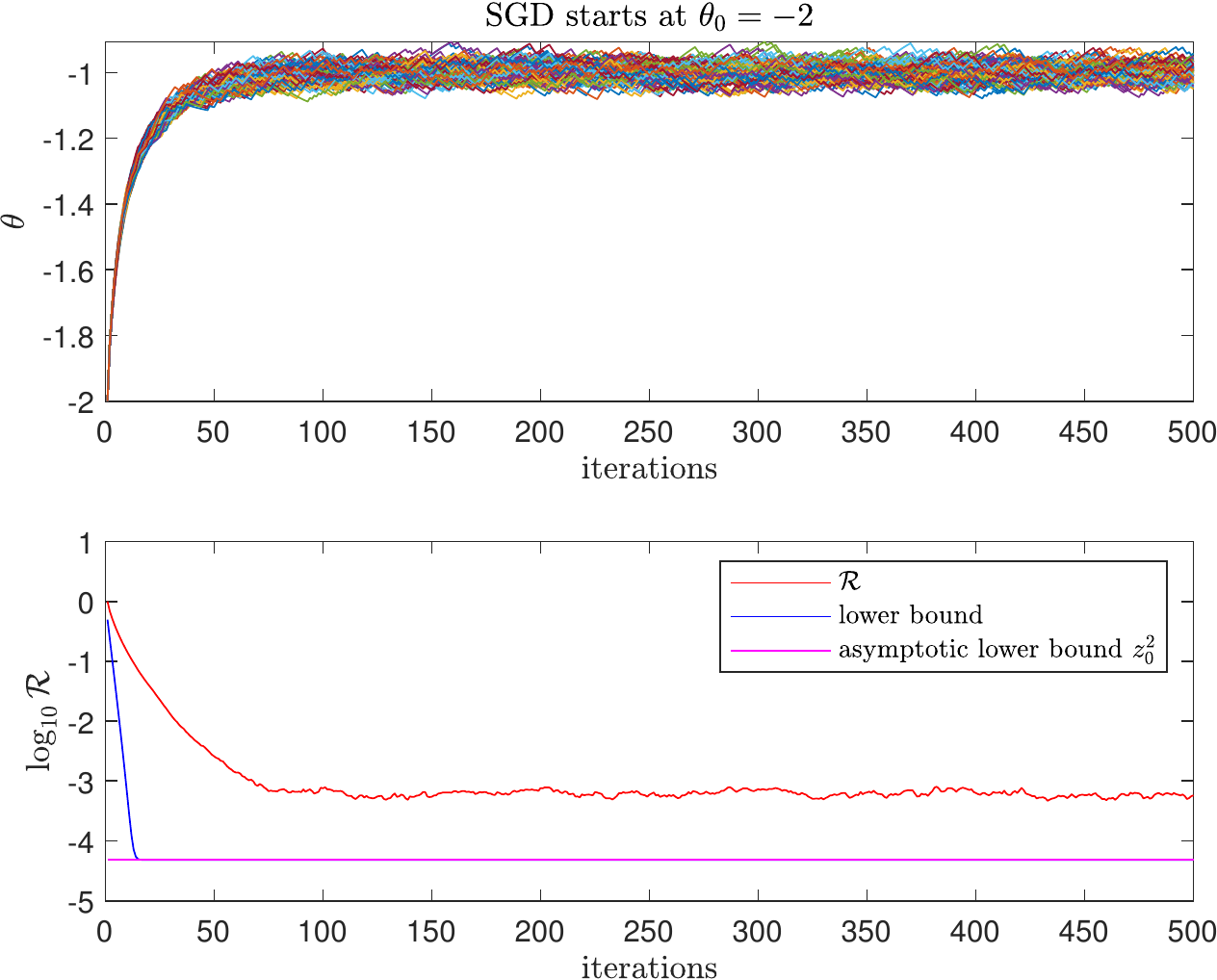}
\caption{$\theta$ trajectories and error curve for SGD starts at $\theta_0 = -2$, $500$ trials.\label{fig:err2_2}}
\end{minipage}
\end{figure}

When running SGD, we set the step-size to be $1/(4\lambda_{\max,0}) \approx 0.069$ and set the maximum iteration to be $500$. The initial point is either $\theta_0 = 2$ or $\theta_0=-2$. For each of the initial points, we run $100$ trials. For each trial, we plot the value of $\theta$ at each iteration.  We show the $\theta$ trajectories of all the $500$ trials and error curves for each initial points in Figures~\ref{fig:err2_1} and~\ref{fig:err2_2}.
In the error curves, the lower bounds are obtained by computing~\eqref{infinequality} in Theorem~\ref{Lemma:iteration}, the asymptotic lower bound $z_0^2$ is the result in Theorem~\ref{ThCB2}.

It is not surprising that SGD will get trapped into local minimum for improper initial point. Figure~\ref{fig:err2_1} shows that when starting at $\theta_0 = 2$, SGD is very likely to get trapped at $\theta = 0.5$. In this case, our lower bound is rather loose. However, as shown in Figure~\ref{fig:err2_2}, even when SGD is trapped at $\theta^\dagger$ the global minimum, there is still an error lower bound, which is very similar to the results in Section~\ref{sec:LR}.

Intuitively, when the objective function is non-convex and has local minima, SGD might get trapped in the local minima. If so, it is obvious that there exists an error lower bound. Our theorems show that even if SGD escapes from the local minimum and approaches the global minimum, there is still a positive error lower bound because of the randomness in constant step-size SGD.

\section{Conclusion}
In this paper, we study the error lower bound of constant step-size SGD. We use expected $L_2$ distance from the global minimum to characterize the error. Our theoretical results show that for potentially non-convex objective function with Lipschitz gradient, there always exists a positive error lower bound for SGD. The lower bounds account for a non-decaying term in classic SGD error upper bound analysis. 



\bibliography{SGD}

\appendix
\appendixpage
\addappheadtotoc

\section{Proof of Theorem 3.2}

\begin{pf}

We first reformulate the right hand side of the iterative formula
\begin{equation}\label{supinequalityapp}
\cR_{k+1}-\cR_k\leq \left[-2\eta\lambda_{\min}+\eta^2\mathsf{\Phi}^2\right]\cR_k+2\eta^2\Lambda D_0\cR^{1/2}_{k}+\eta^2D^2_0\,.
\end{equation}
as a quadratic function, with $\cR_k^{1/2}$ replaced by $x$:
\[
H(x)=\eta\left(\left[-2\lambda_{\min}+\eta\sL^2\right]x^2+2\eta\Lambda D_0x+\eta D^2_0\right)\,.
\]
Since $\eta<\frac{\lambda_{\min}}{\Lambda^2+\lambda^2_{\max,0}}<\frac{\lambda_{\min}}{\Lambda^2+\lambda^2_{\max,0}-\lambda^2_{\min}}$, the coefficient of the quadratic term is negative and the function has a maximum value. By direct calculation, we obtain
\[
\max\,H(x)= \eta^2\left(\frac{\eta\Lambda^2}{2\lambda_{\min}-\eta\sL^2}+D^2_0\right)\leq \left(D^2_0+\frac{1}{2}\right)\eta^2\,,
\]
where we used the condition of $\eta$. Also, the function achieves zero at a root:
\[
x_0=\frac{\eta\Lambda D_0+\sqrt{\eta^2\Lambda^2D^2_0+\eta D^2_0\left(2\lambda_{\min}-\eta\sL^2\right)}}{2\lambda_{\min}-\eta\sL^2}\,.
\]
Therefore, we obtain
\[
H(x_0)\leq \left(D^2_0+\frac{1}{2}\right)\eta^2\,,\quad\text{and}\quad C_0=x_0\,.
\]

Since $H$ is monotonically decreasing in the region of $[x_0\,,\infty)$, one has $H(C_0)\leq 0$. We now use the iterative formula to discuss two cases:
\begin{align*}
\text{if}\quad \cR_k\leq C^2_0:\quad \cR_{k+1}&\leq \cR_{k}+H(W^{1/2}_k)\leq \cR_{k}+\left(D^2_0+\frac{1}{2}\right)\eta^2 \leq C^2_0+\left(D^2_0+\frac{1}{2}\right)\eta^2\,,\\
\text{if}\quad \cR_{k}>C^2_0:\quad \cR_{k+1}&\leq \cR_{k}+H(W^{1/2}_k)\leq \cR_{k}+H(C_0)\leq \cR_{k}\,.
\end{align*}
This implies, for all $k$:
\begin{equation}\label{crKupp}
\cR_k\leq \max\left\{C^2_0,\cR_0\right\}+\left(D^2_0+\frac{1}{2}\right)\eta^2= C_1+\left(D^2_0+\frac{1}{2}\right)\eta^2\,.
\end{equation} 
Plug this into \eqref{supinequalityapp}, we arrive at:
\begin{align*}
\cR_{k+1} \leq & \left[1-2\eta\lambda_{\min}+\eta^2\sL^2\right]\cR_k +2\eta^2\Lambda D_0\left(C_1+\left(D^2_0+1/2\right)\eta^2\right)^{1/2}+\eta^2D^2_0\,,
\end{align*}
and then by setting
\begin{align*}
\alpha & =1-2\eta\lambda_{\min}+\eta^2\sL^2, \\ 
\beta & =2\eta^2\Lambda D_0\left(C_1+\left(D^2_0+1/2\right)\eta^2\right)^{1/2}+\eta^2D^2_0\,,
\end{align*}
we have a simpler version of the updating formula:
\begin{align}
\cR_{k}&\leq \alpha^k\cR_0+\left(\sum^k_{i=0}\alpha^i\right)\beta\leq \alpha^k\cR_0+\frac{\beta}{1-\alpha} =\left(1-2\eta\lambda_{\min}+\eta^2\sL^2\right)^k\cR_0+C_2\eta\,.\label{nonasyminquality}
\end{align}
Similarly, plug \eqref{crKupp} into, 
\begin{equation}\label{infinequalityapp}
\cR_{k+1}-\cR_k\geq \left[-2\eta\lambda_{\max,0}\right]\cR_k-2\eta^2\Lambda D_0\cR^{1/2}_{k}+\eta^2D^2_0\,,
\end{equation}
we also have
\begin{align*}
\cR_{k}\geq & \left[1-2\eta\lambda_{\max,0}\right]\cR_k -2\eta^2\Lambda D_0\left(C_1+\left(D^2_0+1/2\right)\eta^2\right)^{1/2}+\eta^2D^2_0\,,
\end{align*}
and
\[
\cR_{k}\geq \left[1-2\eta\lambda_{\max,0}\right]^k\cR_0+C_3\eta.
\]

Because $\eta<\frac{\lambda_{\min}}{\Lambda^2+\lambda^2_{\max,0}}<\frac{\lambda_{\min}}{\Lambda^2+\lambda^2_{\max,0}-\lambda^2_{\min}}$, we have
\begin{align*}
\alpha &= 1-2\eta\lambda_{\min}+\eta^2\left(\Lambda^2+\lambda^2_{\max,0}-\lambda^2_{\min}\right)\leq 1-\eta\left(2\lambda_{\min}-\eta\left(\Lambda^2+\lambda^2_{\max,0}-\lambda^2_{\min}\right)\right)\leq 1-\eta\lambda_\text{min}\,,
\end{align*}
and setting $K_0=\left\lceil\frac{\log(\eta)}{\log(1-\eta\lambda_{\min})}\right\rceil$ in~\eqref{nonasyminquality}, we have
\[
\begin{aligned}
\cR_{K_0}\leq & \alpha^{K_0}\cR_0+\frac{\beta}{1-\alpha}\\
\leq & \left[1-\eta\lambda_{\min}\right]^{\frac{\log(\eta)}{\log(1-\eta\lambda_{\min})}}\cR_0 +\frac{2\eta\Lambda D_0\left[C_1+\left(D^2_0+1/2\right)\eta^2\right]^{1/2}+\eta D^2_0}{2\lambda_{\min}-\eta\sL^2}\\
= &\left(\cR_0+\frac{2\Lambda D_0\left[C_1+\left(D^2_0+1/2\right)\eta^2\right]^{1/2}+D^2_0}{2\lambda_{\min}-\eta\sL^2}\right)\eta=C_4\eta\,.
\end{aligned}
\]
What's more, for $k\geq K_0$, we also have
\[
\cR_{k}\leq \alpha^{K_0}\cR_0+\frac{\beta}{1-\alpha}\leq C_4\eta\quad  \Rightarrow\quad \cR^{1/2}_{k}\leq C^{1/2}_4\eta^{1/2}\,.
\]

Plug this into the iteration formula \eqref{supinequalityapp} and \eqref{infinequalityapp}, we obtain for $k\geq K_0$
\[
\cR_{k+1}\leq \left[1-2\eta\lambda_{\min}+\eta^2\sL^2\right]\cR_k+2C^{1/2}_4\eta^{5/2}\Lambda D_0+\eta^2D^2_0
\]
\[
\cR_{k+1}\geq \left[1-2\eta\lambda_{\max,0}\right]\cR_k-2C^{1/2}_4\eta^{5/2}\Lambda D_0+\eta^2D_0,
\]
where we also use 
\[
1-2\eta\lambda_{\max,0}>0,\quad \forall \eta<\frac{\lambda_{\min}}{\Lambda^2+\lambda^2_{\max,0}}.
\]
Running the same argument using $\alpha$ and $\beta$ as above, we obtain
\begin{equation}\label{supiteration2app}
\begin{aligned}
\cR_{k} \leq &\left[1-2\eta\lambda_{\min}+\eta^2\mathsf{\Phi}^2\right]^{k-K_0}\cR_{K_0}+\frac{\eta D^2_0+2C^{1/2}_4\eta^{3/2}\Lambda D_0}{2\lambda_{\min}-\eta\mathsf{\Phi}^2}\,,\\
\cR_k \geq &\left[1-2\eta\lambda_{\max,0}\right]^{k-K_0}\cR_{K_0}+\frac{\eta D^2_0-2C^{1/2}_4\eta^{3/2}\Lambda D_0}{2\lambda_{\max,0}}\,.
\end{aligned}
\end{equation}
Finally, the last two inequalites are direct result by letting $k\rightarrow$ in \eqref{supiteration2app}.

\end{pf}

\section{Proof of Theorem 3.3,3.4}

\begin{pf}

First, take $\limsup$ on both sides of \eqref{infinequalityapp}, we can obtain
\begin{align*}
\limsup_{k\rightarrow\infty}\cR_{k+1}\geq & \left[1-2\eta\lambda_{\max,0}\right]\limsup_{k\rightarrow\infty}\cR_k -2\eta^2\Lambda D_0\limsup_{k\rightarrow\infty}\cR^{1/2}_{k}+\eta^2D^2_0\,,
\end{align*}
which implies
\begin{equation}\label{eqn:limsup_nonconv}
2\eta\lambda_{\max,0}\left(\limsup_{k\rightarrow\infty}\cR_k\right)+2\eta^2\Lambda D_0\left(\limsup_{k\rightarrow\infty}\cR_{k}\right)^{1/2}-\eta^2D^2_0\geq 0\,.
\end{equation}
Define a quadratic function to be:
\[
h_1(z)=2\lambda_{\max,0}z^2+2\eta\Lambda D_0z-\eta D^2_0\,.
\]
Comparing it with~\eqref{eqn:limsup_nonconv}, we see that
\[
h_1\left(\left(\limsup_{k\rightarrow\infty}\cR_{k}\right)^{1/2}\right)\geq 0\,.
\]
It is also straightforward to verify that $h_1(z_0)= 0$ by plugging in the definition of $z_0$. Considering $h_1$ is a monotonically increasing function in the region of $[z_0,\infty)$,
\[
\limsup_{k\rightarrow\infty}\cR_k\geq z_0^2\,.
\]

Finally, to prove 
\begin{equation}\label{infb2app}
\liminf_{k\rightarrow\infty}\cR_{k}\geq z^2_0\,.
\end{equation}
and Theorem 3.4, we need to consider two different cases.
\begin{itemize}
\item In the first scenario, we assume there exists $k^*\geq0$ such that
\[
\cR_{k^*}\geq z_0^2\,.
\]
Let
\[
h_2(z)=z^2-h_1(z)=\left[1-2\eta\lambda_{\max,0}\right]z^2-2\eta^2\Lambda D_0z+\eta^2D^2_0\,,
\]
then naturally: $\cR_{k^*+1}\geq h_2(\cR^{1/2}_{k^*})$ according to~\eqref{infinequalityapp}. One can also show that $h_2$ achieves its minimum at 
\[
z^*=\frac{\eta^2\Lambda D_0}{1-2\eta\lambda_{\max,0}}\,,
\] 
then because of condition on $\eta$, $z^\ast<z_0$, and $h_2$ is a monotonically increasing function in the region $[z_0,\infty]$, and thus
\[
\cR_{k^*+1}\geq h_2\left(\cR^{1/2}_{k^*}\right)\geq h_2(z_0)=z_0^2\,.
\]
which implies 
\[
\cR_{k}\geq z^2_0,\quad \forall k>k^*\quad\Rightarrow\quad \liminf_{k\rightarrow\infty}\cR_k\geq z_0^2\,.
\]

\item In the second scenario, we assume for all $k\geq0$, we have
\[
0\leq \cR_{k}\leq z_0^2\quad\Rightarrow\quad h_1\left(\cR^{1/2}_{k}\right)\leq h_1(z_0)=0\,.
\]
Here we used the fact that $h_1$ is monotonically increasing in the region of $[0,z_0]$. Using~\eqref{infinequalityapp}, we have, for all $k\geq 0$:
\[
z_0^2\geq \cR_{k+1}\geq \cR_{k}-h_1\left(\cR^{1/2}_{k}\right)\geq \cR_k-h_1(z_0)=\cR_{k}\,,
\]
meaning $\{\cR_k\}$ is an increasing sequence with an upper bound $z^2_0$. However, we also have, according to  $\limsup_{k\rightarrow\infty}\cR_k\geq z^2_0$, therefore, we finally obtain
\[
\lim_{k\rightarrow\infty}\cR_k=z_0^2\,.
\]
\end{itemize}
Combining the discussion of the two scenarios, we prove \eqref{infb2app} and Theorem 3.4.
\end{pf}

\end{document}